Original Article

# Transitive perfect colorings of the non-regular Archimedean tilings

**Junmar Gentuya[1],* and René Felix[2]**

[1] *Division of Natural Sciences and Mathematics*
*University of the Philippines Visayas Tacloban College*
*Magsaysay Blvd, Tacloban City*
*Leyte Province 6500, Philippines*

[2] *Institute of Mathematics*
*College of Science*
*University of the Philippines Diliman*
*Quezon City 1101, Philippines*

**ABSTRACT.** In this work, we give a method to obtain nontrivial transitive perfect colorings of the non-regular Archimedean tilings using the least possible number $n$ of colors. We also look for other non-equivalent transitive perfect $n$-colorings of a given non-regular Archimedean tiling.

## Introduction

In 1993, R. L. Roth considered the problem of obtaining nontrivial perfect colorings of multipatterns in the Euclidean plane. A *multipattern* is a symmetrical structure in the Euclidean plane, with symmetry group $G$, which has several $G$-orbits of motifs. Moreover, Roth assumed that $G$ acts transitively on the colors. Given a multipattern, the questions he answered in his paper are the following: (1) can we obtain a nontrivial perfect coloring of the given multipattern such that $G$ acts transitively on the set of colors? and (2) what is the minimum number of colors that can be used to obtain nontrivial perfect colorings of a given multipattern?

In Table 1, Roth's coloring numbers of all 17 types of plane crystallographic groups are presented. If $G$ is the symmetry group of a multipattern, then $N(G)$, called the *coloring number* of $G$, is the minimum number (> 1) of colors which suffice to color the multipattern. Roth considered a "worst-case scenario" of a given multipattern with symmetry group $G$. This means that the motifs of a multipattern are placed using all possible conjugacy classes of stabilizer subgroups of $G$. So, regardless of how a multipattern looks like, as long as its symmetry group

*Corresponding author:

*Division of Natural Sciences and Mathematics*
*University of the Philippines Visayas Tacloban College*
*Magsaysay Blvd, Tacloban City*
*Leyte Province 6500, Philippines*
*Email Address: jggentuya@up.edu.ph*

Table 1. The plane crystallographic groups and their coloring numbers (Roth 1993)

| Symmetry Group $G$ | Coloring Number $N(G)$ |
|---|---|
| $p1$ | 1 |
| $pg$ | 2 |
| $pm$ | 2 |
| $cm$ | 2 |
| $pgg$ | 2 |
| $p2$ | 3 |
| $pmg$ | 3 |
| $pmm$ | 3 |
| $cmm$ | 3 |
| $p3$ | 4 |
| $p3m1$ | 4 |
| $p31m$ | 4 |
| $p4$ | 5 |
| $p4m$ | 9 |
| $p4g$ | 9 |
| $p6$ | 7 |
| $p6m$ | 25 |

$G$ is a plane crystallographic group, the multipattern can always be perfectly colored using $N(G)$ colors.

## Tilings Under Consideration

The list of the 11 Archimedean tilings is shown in Table 2. Three of these tilings are called *regular* tilings (see Fig. 1). They are tilings by equilateral triangles, squares and regular hexagons. The tiles of each of the regular tilings form only one $G$-orbit of tiles. This means that regular tilings are not multipatterns. Grünbaum and Shephard (1977) studied the perfect colorings of regular tilings. The non-regular Archimedean tilings (see Fig. 2) are examples of multipatterns where each is made up of tiles congruent to 2 or 3 regular polygons. The number of G-orbits of tiles of these tilings can be seen in Table 2. Santos and Felix (2006) studied the perfect colorings of a given non-regular Archimedean tiling with 2 orbits of colors.

Table 2. Archimedean tilings: their respective symmetry groups and numbers of $G$-orbits of tiles

| Archimedean Tiling | Symmetry Group | Number of $G$-orbits |
|---|---|---|
| $3 \cdot 3 \cdot 3 \cdot 3 \cdot 3 \cdot 3$ | $p6m$ | 1 |
| $4 \cdot 4 \cdot 4 \cdot 4$ | $p4m$ | 1 |
| $6 \cdot 6 \cdot 6$ | $p6m$ | 1 |
| $3 \cdot 3 \cdot 3 \cdot 4 \cdot 4$ | $cmm$ | 2 |
| $3 \cdot 3 \cdot 3 \cdot 3 \cdot 6$ | $p6$ | 3 |
| $3 \cdot 3 \cdot 4 \cdot 3 \cdot 4$ | $p4g$ | 2 |
| $4 \cdot 8 \cdot 8$ | $p4m$ | 2 |
| $3 \cdot 6 \cdot 3 \cdot 6$ | $p6m$ | 2 |
| $3 \cdot 12 \cdot 12$ | $p6m$ | 2 |
| $3 \cdot 4 \cdot 6 \cdot 4$ | $p6m$ | 3 |
| $4 \cdot 6 \cdot 12$ | $p6m$ | 3 |

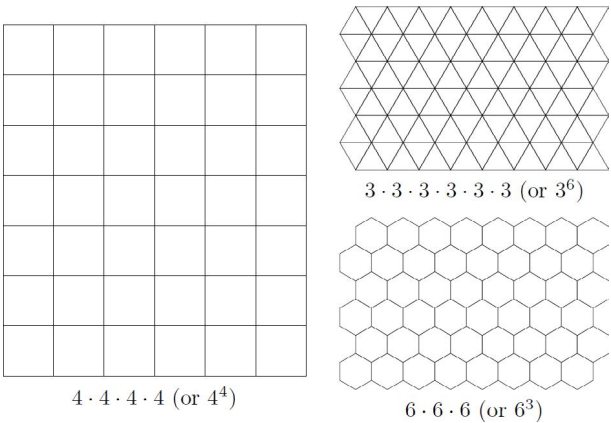

Figure 1. The regular Archimedean tilings.

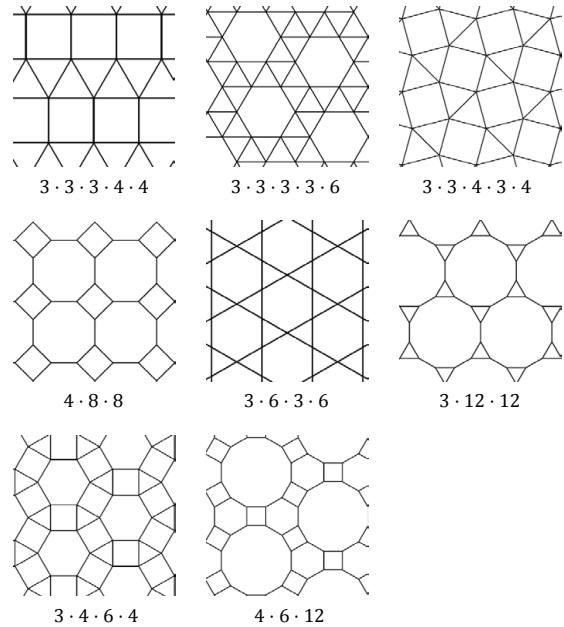

Figure 2. The non-regular Archimedean tilings.

**Objectives**

Our aim in this work is to look for nontrivial transitive perfect colorings of each non-regular Archimedean tiling using the least possible number of colors (the coloring number of the tiling). We shall give a method based on a coloring framework (Felix 2011) that would guide us to obtain a nontrivial transitive perfect coloring of a given non-regular Archimedean tiling.

**Preliminaries**

Let $\mathbb{X}$ be a plane. An *isometry* $f$ of $\mathbb{X}$ is a bijective map $f : \mathbb{X} \to \mathbb{X}$ that preserves distance. Let $X \subset \mathbb{X}$ ($X$ can be any set of points or a tiling in $\mathbb{X}$). If $f(X) = X$, then $f$ is called a *symmetry* of $X$. All symmetries of $X$ form a group under a composition called the *symmetry group* of $X$, denoted usually by $G$.

Let $G$ be a group acting on a set $X$. If $x \in X$, then $Gx \coloneqq \{gx : g \in G\}$ is called the $G$-orbit of $x$ (this set contains the images of $x$ under the elements of $G$), and $Stab_G(x) \coloneqq \{g \in G : gx = x\}$ is called the *stabilizer* of $x$ in $G$.

**Proposition 1.** Let $G$ be a group acting on a set $X$. If $x_1, x_2 \in X$ such that $gx_1 = x_2$ for some $g \in G$, then $Stab_G(x_2) = g(Stab_G(x_1))g^{-1}$.

*Proof.* Note that
$$\begin{aligned}
h \in Stab_G(x_2) &\Leftrightarrow hx_2 = x_2 \\
&\Leftrightarrow h(gx_1) = gx_1 \\
&\Leftrightarrow g^{-1}hgx_1 = x_1 \\
&\Leftrightarrow g^{-1}hg \in Stab_G(x_1) \\
&\Leftrightarrow h \in g(Stab_G(x_1))g^{-1}
\end{aligned}$$
Therefore, $Stab_G(x_2) = g(Stab_G(x_1))g^{-1}$. ∎

Under the action of a symmetry group $G$, in general, a set is decomposed into several subsets called $G$-orbits or *transitivity classes*. Proposition 1 tells us that if $x_1$ and $x_2$ belong to one $G$-orbit, then their stabilizers in $G$ are conjugate. If the $G$-orbit of $x \in X$ is the entire set $X$, then we say that $G$ acts transitively on $X$. On the other hand, if $X$ is a union of at least 2 distinct $G$-orbits of elements of $X$, then $X$ is a *multipattern*.

**Archimedean Tilings**

A tiling of the plane $\{t_1, t_2, ..., t_n, ...\}$ is a countable family of connected closed subsets of the plane called *tiles*, such that the union of the sets $t_i$ is the entire plane and the interiors of distinct sets $t_i$ and $t_j$ are disjoint. The point of intersection in a given tiling where 3 or more tiles meet is called a *vertex*. The collection of points shared by 2 tiles is called an *edge*. If each side of every tile is an edge of the tiling and vice versa, then the tiling is called *edge-to-edge*.

The Archimedean $p_1 \cdot p_2 \cdot ... \cdot p_q$ tiling is an edge-to-edge tiling composed of tiles which are regular polygons, with a $p_1$-gon, a $p_2$-gon, ..., and a $p_q$-gon arranged in a cyclic manner around each and every vertex such that the symmetry group of the tiling acts transitively on its vertices.

### Colorings

Let $X$ be a set of objects to be colored and let $G$ be its symmetry group. If $C = \{c_1, c_2, ..., c_n\}$ is a set of $n$ colors, then an $n$-coloring of $X$ is a surjective (onto) map $c : X \rightarrow C$. Each object $x$ in $X$ is assigned a color $c(x) \in C$. The $n$-coloring of $X$ is given by the partition $P = \{P_1, P_2, ..., P_n\}$ of $X$ where 2 elements of $X$ belong to same set in $P$, if and only if they are assigned the same color.

An element in $G$ that effects permutation of colors is called a *color symmetry*. Color symmetries form a subgroup of $G$ called the color symmetry group or simply, *color group*, denoted by $H$. A coloring of $X$ is said to be perfect if $G = H$, i.e., every element of $G$ is a color symmetry. A perfect coloring is said to be *trivial* if only one color is used in the coloring. Meanwhile, a coloring is *transitive* if $G$ acts transitively on the colors, that is, if $c_i$ and $c_j$ are colors in $C$, then there is an element of $G$ sending $c_i$ to $c_j$. In a transitive coloring of $X$, the set of all objects in $X$ of one color is congruent to another set of all objects in $X$ of another color.

Two colorings of the same symmetrical structure are called equivalent if one of the colorings is obtained from the other by (1) a bijection from $C_1$ (set of colors used in the first coloring) to $C_2$ (set of colors used in the second coloring), or (2) a symmetry in the symmetry group $G$ of the uncolored symmetrical structure, or (3) a combination of (1) and (2).

### Coloring Framework

Let $X$ be a set and $G$ a group acting transitively on $X$. The Coloring Framework (Felix 2011) gives a method for finding all possible colorings of a set $X$ for which the group $G$ permutes the colors of $C$. We state this as follows (see Evidente 2012 for the proof).

**Theorem 1**. Let $X$ be a set and let $G$ be a group acting transitively on $X$.

1. If $\{P_i\}_{i=1}^n$ is a coloring of $X$ for which $G$ permutes the colors, then for every $x \in X$, there exists $J \leq G$ such that $Stab_G(x) \leq J$, and the coloring is described by the partition $\{gJx : g \in G\}$.

2. Let $x \in X$ and $J \leq G$ such that $Stab_G(x) \leq J$ and $[G : J] = n < \infty$. Then $P = \{gJx : g \in G\}$ is an $n$-coloring of $X$ for which $G$ permutes the colors.

As a remark, it is assumed in the theorem that group $G$ acts transitively on set $X$. Hence, the resulting n-coloring of $X$ given by the partition $P$ is perfect and transitive.

Note that the symmetry group $G$ of a given Archimedean tiling acts transitively on each $G$-orbit of tiles in the tiling. Based on the theorem, we can obtain an $n$-coloring of a $G$-orbit which is perfect and transitive using $J$ such that $Stab_G(x) \leq J \leq G$ where $[G : J] = n$, and $x$ is in this $G$-orbit. Each $G$-orbit will be colored independently. Then the colored $G$-orbits will be combined but we intend to have a transitively and perfectly colored Archimedean tiling as a result.

In the work of Santos and Felix (2006), a perfect coloring of a given Archimedean tiling is obtained by doing the following: (1) color perfectly each $G$-orbit of tiles using a subgroup of $G$, and (2) combine the colored $G$-orbits of tiles. To perfectly color a $G$-orbit $X_i$, a subgroup $J_i$ is used. The condition needed to be satisfied is that $J_i$ must contain the stabilizer of a tile in $X_i$. To color perfectly another $G$-orbit $X_k$ of tiles, another subgroup $J_k$ is used, that is, a subgroup that contains the stabilizer of a tile in $X_k$. The subgroup $J_k$ can be different from $J_i$. In effect, in the resulting perfect coloring of the tiling, the $G$-orbits of tiles do not share colors. Thus, there are 2 orbits of colors used: one for $X_i$ and another for $X_k$.

In this work, an interest is to have the symmetry group $G$ to be transitive on the colors. In other words, only one set of colors is used to color perfectly each and every $G$-orbit of tiles. This means that the $G$-orbits of tiles share colors. The idea is as follows: To color perfectly a $G$-orbit $X_i$, we use a subgroup $J$ in $G$ if it contains $Stab_G(t_i)$ where $t_i \in X_i$. If a color $c$ is used to color the tile $t_i$, then $c$ is also used to color all the tiles belonging to $Jt_i$, the $J$-orbit of $t_i$. We shall use $J$ in coloring perfectly $X_k$, another $G$-orbit of tiles, if $J$ contains the stabilizer of a tile $t_k \in X_k$ and we will also use the same color $c$. Thus, color $c$ is assigned to the tiles in $Jt_k$. In other words, $c$ is assigned to tiles in the set $J\{t_i, t_k\} = Jt_i \cup Jt_k$.

Let $G$ be the symmetry group of an Archimedean tiling. Note that the list of subgroups of $G$ up to conjugacy can be obtained using the computer algebra system GAP (2008). We consider Roth's coloring number of $G$ as the upper bound for the number of colors to be used in coloring perfectly and transitively the given tiling. The focus of the investigation is finding a subgroup $J$ of $G$ whose index $n$ in $G$ is possibly lower than Roth's coloring number of $G$. We look for other non-equivalent transitive perfect $n$-colorings of the given Archimedean tiling.

Below we outline the steps to obtain a nontrivial transitive perfect coloring of an Archimedean tiling using the least possible number of colors.

1. Given an Archimedean tiling, determine the symmetry group $G$.

2. Determine the $G$-orbits of tiles in the tiling.

3. Suppose that there are m number of $G$-orbits of tiles. Pick a tile $x_i$ in the $G$-orbit $X_i$ of tiles, for each $i=1,2,...,m$.

4. Get $Stab_G(x_i)$ for each $i = 1,2,...,m$. Note that this subgroup of $G$ is isomorphic to a finite group $C_k$ (cyclic group of order $k$) or $D_k$ (dihedral group of order $2k$) for some integer $k$.

5. Find a proper subgroup $J$ of smallest possible index $n$ in $G$ such that for every $i = 1,2,...,m$, a conjugate of $Stab_G(x_i)$ is contained in $J$.

6. Determine the $J$-orbits of tiles in $X_i$, $i = 1,2,...,m$.

7. For each $i = 1,2,...,m$, pick $t_i$ in $X_i$ such that $Stab_G(t_i) \leq J$. Form the set $T = \{t_1, t_2, ..., t_m\}$ (the set of all $J$-orbit representatives $t_i$ from $X_i$, $i = 1,2,...,m$).

8. Let $\{g_1, g_2, ..., g_n\}$ be a complete set of left coset representatives of $J$ in $G$. Form the partition $P = \{gJT : g \in G\} = \{g_1JT, ..., g_nJT\}$.

9. Let $C = \{c_1, c_2, ..., c_n\}$ be a set of colors. To each set $g_jJT$ in $P$, assign the color $c_j$, $j = 1,2,...,n$.

The partition $P$ corresponds to a nontrivial transitive perfect $n$-coloring of the given Archimedean tiling.

## Results

We shall apply the method to the tiling $3 \cdot 3 \cdot 3 \cdot 4 \cdot 4$ and obtain its transitive perfect colorings using the least number of colors. For the rest of the non-regular Archimedean tilings, we shall omit the details but the resulting transitive perfect colorings will be summarized in a table.

Tiling $3 \cdot 3 \cdot 3 \cdot 4 \cdot 4$ is shown in Figure 2. Around each vertex of the tiling are 3 equilateral triangles and 2 squares. Consider the square tile colored yellow in Figure 3 and call it $t$. The symmetry group of this tiling is $G = \langle u, v, r, s \rangle \cong cmm$, where $u$ and $v$ are 2 linearly independent translations described by the vectors indicated in the figure, $r$ is a mirror reflection about a horizontal line passing through the center of $t$, and $s$ is a mirror reflection about a vertical line passing through the center of $t$. Under the action of the symmetry group $G$, the squares form one $G$-orbit $X_1$, and the triangles form another $G$-orbit $X_2$. The stabilizer in $G$ of a square tile is a group of type $D_2$, a dihedral group of order 4. For example, the tile $t$ is stabilized by the group $\langle r, s \rangle \cong D_2$. On the other hand, the stabilizer in $G$ of a triangle is of type $D_1$, a dihedral group of order 2.

The possible subgroups (of finite index) of a plane crystallographic group of type $cmm$ are of types $p1, p2, pm, pg, cm, pmm, pmg, pgg,$ and $cmm$ (Rapanut 1988). Among them, only plane crystallographic groups of type $pmm$ or $cmm$ could contain subgroups of types $D_1$ and $D_2$.

The coloring number of $cmm$ is 3 (Roth 1993). Table 3 lists the subgroups of $cmm$ of index less than or equal to 3 where each has $r$ and $s$ among the generators. It is easy to see that each of these subgroups is either of type $pmm$ or of type $cmm$.

A subgroup of type $pmm$ in the symmetry group $G \cong cmm$ has index $2n$ (Rapanut 1988), where $n \in \mathbb{Z}$. When $n = 1$, we have a subgroup isomorphic to $pmm$ of index 2 in $G$ and this is the subgroup $\langle u^2, uv, r, s \rangle$. Note that in Step 5 of the method, we need a proper subgroup of smallest possible index in $G$ that contains a conjugate of a stabilizer of a tile in each $G$-orbit. The subgroup $\langle u^2, uv, r, s \rangle$ has such property: it contains a conjugate of the stabilizer of a square in $X_1$ and of a triangle in $X_2$. Hence, this will be the subgroup $J$.

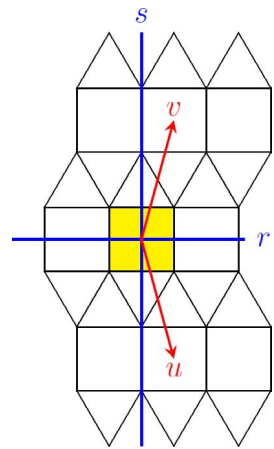

Figure 3. Tiling $3 \cdot 3 \cdot 3 \cdot 4 \cdot 4$ with generators of $G$.

Table 3. Subgroups of type pmm or cmm of index at most 3 in $G = \langle u, v, r, s \rangle \cong cmm$

| Subgroup of $G$ | Index in $G$ | Type of Group |
|---|---|---|
| $\langle u, v, r, s \rangle$ | 1 | cmm |
| $\langle u^2, uv, r, s \rangle$ | 2 | pmm |
| $\langle u^{-1}v, r, s, u^3 \rangle$ | 3 | cmm |
| $\langle uv, r, s, u^3 \rangle$ | 3 | cmm |

We then take the $J$-orbits of squares in $X_1$ and triangles in $X_2$. We superimpose the reduced diagram of the subgroup $J$ on the tiling as shown in Figure 4a, so that the $J$-orbits of tiles are easily obtained as shown in Figures 4b and 4c. Moreover, observe that any square in $X_1$ has stabilizer in $G$ contained in $J$ and any triangle in $X_2$ has stabilizer in $G$ contained in $J$. Since $[G:J] = 2$, the least of all indices of proper subgroups of $G$, $J$ induces a coloring of the tiling that uses 2 colors. Thus, the coloring number of the $3 \cdot 3 \cdot 3 \cdot 4 \cdot 4$ tiling is 2.

Let $t_1$ be any yellow or gray square in Figure 4b and $t_2$ be any cyan or magenta triangle in Figure 4c. Then form the set $T = \{t_1, t_2\}$, which consists of $J$-orbit representatives from $X_1$ and $X_2$, respectively, such that $Stab_G(t_1), Stab_G(t_2) \leq J$. Table 4 gives the summary of the possible combinations of $J$-orbit representatives from $X_1$ and $X_2$. With the set $\{e, u\}$ as a complete set of left coset representatives of $J$ in $G$, we form the partition $\{gJT : g \in G\} = \{gJ\{t_1, t_2\} : g \in G\} = \{J\{t_1, t_2\}, uJ\{t_1, t_2\}\}$. We assign color red to $J\{t_1, t_2\}$ and green to $uJ\{t_1, t_2\}$. For example, in Table 4 row (a), if $t_1$ is a yellow square in $X_1$ and $t_2$ is a cyan triangle in $X_2$, then the resulting coloring is shown in Figure 5a. This is a transitive perfect 2-coloring of the tiling. However, using $J$-orbit representatives in Table 4 row (d), the resulting coloring is equivalent to Figure 5a. Moreover, the transitive perfect 2-colorings (see Figure 5b) obtained by using $J$-orbit representatives in Table 4 row (b) and row (c) are equivalent. Therefore, there are 2 non-equivalent transitive perfect 2-colorings of the tiling $3 \cdot 3 \cdot 3 \cdot 4 \cdot 4$.

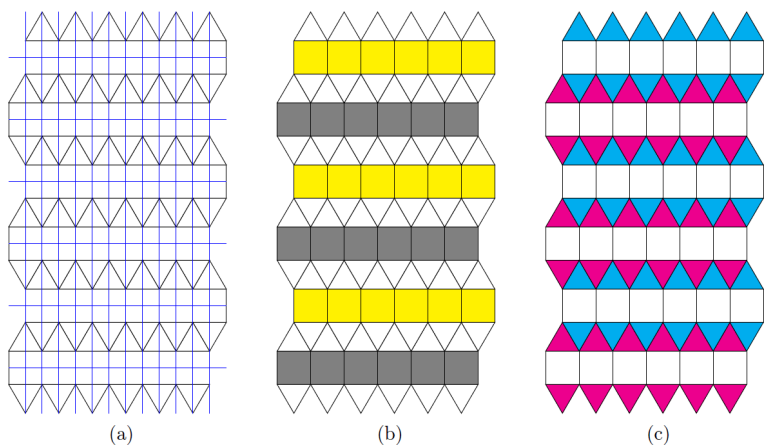

Figure 4. Tiling $3 \cdot 3 \cdot 3 \cdot 4 \cdot 4$ overlaid by the reduced diagram of $J$; The $J$-orbits of (b) squares in $X_1$ and (c) triangles in $X_2$.

Table 4. Combinations of J-orbit representatives of tiles in $3 \cdot 3 \cdot 3 \cdot 4 \cdot 4$.

| Coloring | Representative from $X_1$ | Representative from $X_2$ |
|---|---|---|
| (a) | yellow | cyan |
| (b) | yellow | magenta |
| (c) | gray | cyan |
| (d) | gray | magenta |

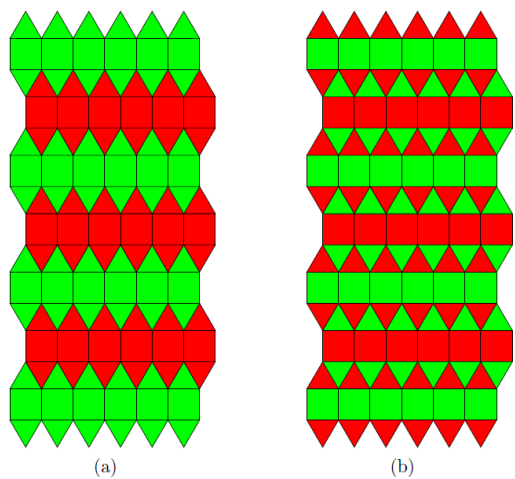

Figure 5. The 2 transitive perfect 2-colorings of the tiling $3 \cdot 3 \cdot 3 \cdot 4 \cdot 4$.

For the rest of the non-regular Archimedean tilings, we also applied the method to find nontrivial transitive perfect colorings (see Table 5). These colorings are shown in Figures 6–8. Aside from $3 \cdot 3 \cdot 3 \cdot 4 \cdot 4$ tiling whose coloring number 2 does coincide with Roth's result, we also have the $3 \cdot 3 \cdot 3 \cdot 3 \cdot 6$ tiling with symmetry group of type *p6* whose coloring number is 4 and not 7 (Roth's coloring number of *p6*). As for the $3 \cdot 6 \cdot 3 \cdot 6$ and $3 \cdot 12 \cdot 12$ tilings, both with symmetry group of type *p6m*, their coloring number is 4 and not 25 (Roth's coloring number of *p6m*). Note that these tilings are *not* multipatterns in the worst-case scenario described in Roth's article. These tilings do not have all types of stabilizer subgroups. One way to look at these is through the reduced diagrams of the tilings. Take for instance the $3 \cdot 3 \cdot 3 \cdot 4 \cdot 4$ tiling. Its reduced diagram consists of (a) intersecting axes of mirror reflections where two (non-conjugate) reflections have a product which is a half-turn, and (b) centers of half-turns not lying on the axes of mirror reflections. A subgroup of *cmm* that fixes a center of half-turns in (b) above is of type $C_2$ (cyclic group of order 2); however, no tile in the tiling is stabilized by $C_2$; Otherwise, the result coincides with that of Roth's.

## Conclusion and Recommendations

The method that we provide in this paper was applied to find nontrivial inequivalent transitive perfect colorings of the non-regular Archimedean tilings. The inequivalent transitive perfect colorings of these tilings are also presented. One may consider finding possibly other transitive perfect colorings of Archimedean tilings using a number of colors more than the coloring numbers obtained in this work or even more than Roth's coloring numbers. We suggest developing a computer program that can generate these colorings. This may lead one to construct a method that systematically counts inequivalent transitive perfect colorings of Archimedean tilings possibly without doing the actual colorings or looking at different combinations of orbits of tiles. One might find insights by studying first the work of Frettlöh (2008), and that of Santos and Felix (2006) and the references therein.

Transitive perfect colorings of $k$-uniform tilings in the Euclidean plane are also interesting to investigate. As a matter of fact, transitive perfect colorings of 2-uniform tilings (in the Euclidean plane) were already done by Felix and Eclarin (2014) recently, so one can consider $k$-uniform tilings where $k \geq 3$.

## Acknowledgement

The first author gratefully acknowledges the financial assistance given by DOST-SEI through the Accelerated Science and Technology Human Resource Development Program (ASTHRDP) and University of the Philippines Visayas.

Table 5. The least number of colors used and the number of inequivalent transitive perfect colorings of the 7 other non-regular Archimedean tilings

| Archimedean Tiling | Least number $n$ of Colors | Number of Inequivalent $n$-Colorings |
|---|---|---|
| $3 \cdot 3 \cdot 3 \cdot 3 \cdot 6$ | 4 | 4 |
| $3 \cdot 3 \cdot 4 \cdot 3 \cdot 4$ | 9 | 3 |
| $4 \cdot 8 \cdot 8$ | 9 | 1 |
| $3 \cdot 6 \cdot 3 \cdot 6$ | 4 | 1 |
| $3 \cdot 12 \cdot 12$ | 4 | 1 |
| $3 \cdot 4 \cdot 6 \cdot 4$ | 25 | 1 |
| $4 \cdot 6 \cdot 12$ | 25 | 1 |

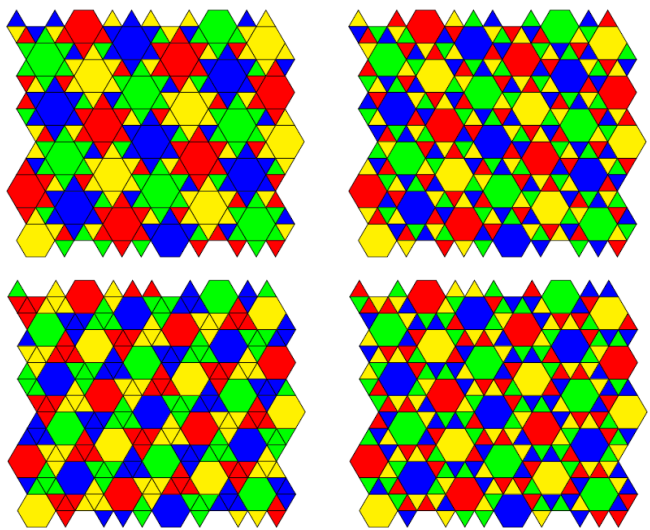

Figure 6. The 4 transitive perfect 4-colorings of the tiling $3 \cdot 3 \cdot 3 \cdot 3 \cdot 6$.

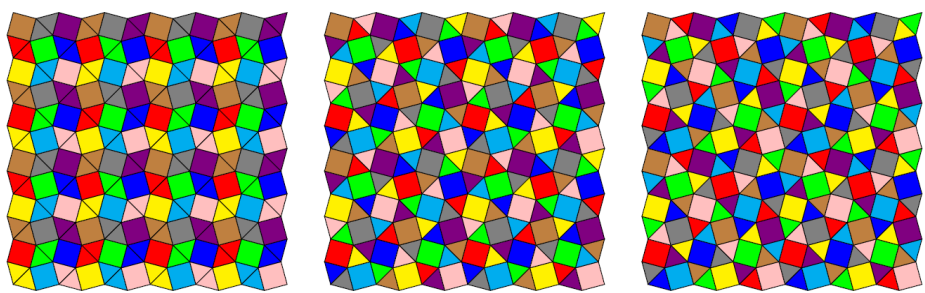

Figure 7. The 3 transitive perfect 9-colorings of $3 \cdot 3 \cdot 4 \cdot 3 \cdot 4$.

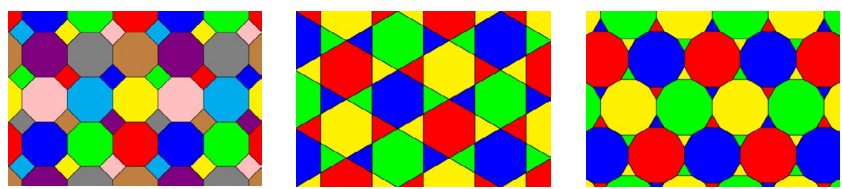

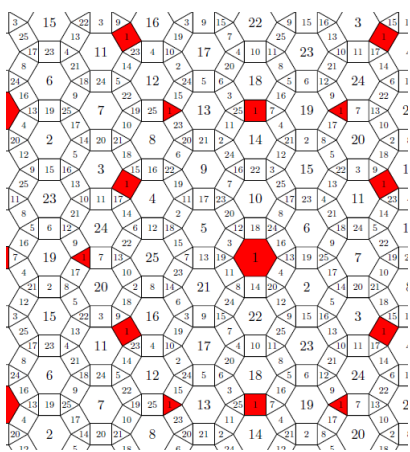
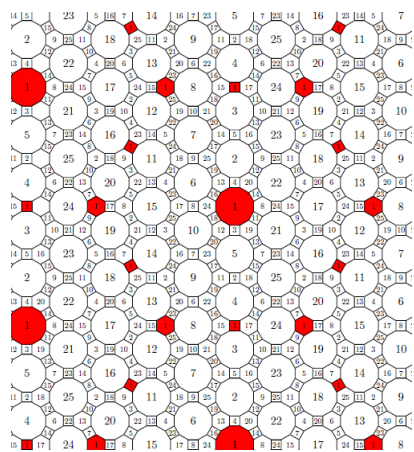

Figure 8. The transitive perfect (a) 9-coloring of $4 \cdot 8 \cdot 8$, (b) 4-coloring of $3 \cdot 6 \cdot 3 \cdot 6$, (c) 4-coloring of $3 \cdot 12 \cdot 12$, (d) 25-coloring of $3 \cdot 4 \cdot 6 \cdot 4$, (e) 25-coloring of $4 \cdot 6 \cdot 12$.